\documentclass[12pt]{amsart}
\usepackage[english]{babel}
\usepackage[all]{xy}
 \usepackage[latin1]{inputenc}
  \usepackage[T1]{fontenc}
  \usepackage{amsmath,amssymb}
  \usepackage{amsmath}
  \usepackage{amsthm}
\usepackage{amssymb}
\usepackage{mathrsfs}
\usepackage{enumerate}
\usepackage{graphicx}
\usepackage{tikz}
\usetikzlibrary{shapes.geometric}
\usetikzlibrary{positioning}
\usetikzlibrary{arrows}
\usetikzlibrary{shapes.multipart}
\usepackage{physics}
\usepackage[notcite, final, notref]{showkeys}
\usepackage{dsfont}
\DeclareMathOperator*{\esssup}{ess\,sup}

\newtheorem{theorem}{Theorem}[section]

\newtheorem{lemma}[theorem]{Lemma}
\newtheorem{proposition}[theorem]{Proposition}

\newtheorem{corollary}[theorem]{Corollary}
\newtheorem{definition}[theorem]{Definition}

\theoremstyle{definition}
\newtheorem{remark}[theorem]{Remark}

\newtheorem{example}[theorem]{Example}

\newcommand{\e}{\varepsilon}

\usepackage{tikz}
\usetikzlibrary{arrows,shapes,positioning,shadows,trees,matrix}
\usepackage{xcolor}

\begin{document}

\title[Betweenness of membership functions]{Betweenness of membership functions: classical case and hyperbolic-valued functions}

\author[D. Alpay]{Daniel Alpay}
\address{(DA) Schmid College of Science and Technology \\
Chapman University\\
One University Drive
Orange, California 92866\\
USA}
\email{alpay@chapman.edu}

\author[L. Mayats-Alpay]{Liora Mayats-Alpay}
\address{(LMA) Fowler School of Engineering and Schmid College of Science and Technology\\
Chapman University\\
One University Drive
Orange, California 92866\\
USA}
\email{mayatsalpay@chapman.edu}

\thanks{D. Alpay thanks the Foster G. and Mary McGaw Professorship in
  Mathematical Sciences, which supported his research.}

\maketitle

\begin{abstract}
We study betweenness of membership functions in the fuzzy setting and for membership functions taking values in the set of hyperbolic numbers.
\end{abstract}

\tableofcontents
\section{Introduction}
\setcounter{equation}{0}

{\bf Prologue:}
To develop mathematical tools to study similarity of objects or situations is a very important problem in a wide range of topics, from botany to
psychology and more, and involve in particular comparisons of sets, finite or infinite.
We mention for instance Paul Jaccard \cite{jaccard1901distribution}, whose studies
of comparative floral distribution lead to the notion of Jaccard index of similarity ({\sl coefficient de communaut\'e}, in French).
To define this index,  some notations need to be introduced. Given a set  $\Omega$, we denote by $\Omega\setminus A$ the complement of
$A\in\mathcal P(\Omega)$ in $\Omega$ and by $A\Delta B$ the symmetric difference of $A$ and $B$:
\begin{equation}
A\Delta B=(A\setminus B)\cup(B\setminus A)=(A\cup B)\setminus (A\cap B).
\end{equation}
The Jaccard index is defined by
\begin{equation}
J(A,B)=\frac{{\rm Card}\, A\cap B}{{\rm Card}\, A\cup B}
  \end{equation}
  where $A$ and $B$ are two finite subsets (not both empty) of a common set $\Omega$. If  $A=B=\emptyset$, one sets $J$ to be $0$.
  Note that $J(A,B)\in[0,1]$ and that
  \begin{equation}
    \label{dist-jac}
    D(A,B)=1-J(A,B)=\frac{{\rm Card}\, A\cup B-{\rm Card}\, A\cap B}{{\rm Card}\, A\cup B}
    =\frac{{\rm Card}\, A\Delta B}{{\rm Card}\, A\cup B}
    \end{equation}
    is a well known distance on sets,\footnote{See Definition \ref{metric-space} for the notion of distance.} similar to the distance between finite random variables introduced earlier in information theory by
    C. Rajski; see e.g. \cite{MR0345715} and \cite{rajski_1961} for the
    latter. The proofs in these papers are easily adapted to the case of finite sets.
    We also note the works \cite{marczewski1958certain} and \cite{zbMATH03194216}\footnote{in Polish; available online.} of Marczewski and Steinhaus, where \eqref{dist-jac} (and counterparts
    for functions)
    is investigated, with applications to study of species growing in forests (analysis of biotopes). Note that these authors do not mention Jaccard.
    For another later proof of the triangle inequality for the Jaccard distance, see e.g. \cite{levandowsky1971distance}.\smallskip

    We also mention the work of Amos Tversky \cite{tversky1977features}, where  a representation theorem to measure similarity between different sets of objects is developed
using the decomposition of $A\cup B$ into the non-overlapping sets $A\setminus B$, $B\setminus A$ and $A\Delta B$. For a
recent survey we refer to \cite{MR2522441}.\smallskip

{\bf Betweenness:}
In the study of similarity of sets, estimating the betweenness of a set of features with respect to two other sets of features is a major question, which can be defined and studied in
different ways, depending on the underlying
structure. It involves important analytic tools, such as metric spaces, strictly convex norms and lattices. In a general metric space one can define the
notion as corresponding to cases of equality in the triangle inequality. In a vector space it is easy to define betweenness:  a vector is between two
vectors $u$ and $v$ if it belongs to the closed interval  defined by these
two vectors (or equal to $u$ when $u=v$).  In \cite{restle1959metric} Restle defines and studies the notion of betweenness of sets.
Let $\Omega$ be a set and
let $A,B,C\in\mathcal P(\Omega)$. Following Restle (see \cite[Definition 2 p. 210]{restle1959metric}) one says that the set $C$ is between the sets $A$
and $B$ if\footnote{Restle writes these two inclusion conditions in a slightly different, but equivalent, way.}
\begin{equation}
\label{crispbetween}
A\cap B\subset C\subset A\cup B,
\end{equation}
which can be translated in terms of indicator functions (see \eqref{indic}) as
\begin{equation}
\label{ineq456}
  1_{A\cap B}(x)\le 1_C(x)\le 1_{A\cup B}(x).
  \end{equation}
Among other questions Restle is interested in \cite{restle1959metric} in the case of equality in the triangle inequality in an
underlying metric space.  To palliate the lack of vector space structure Restle introduces the notion of linear array of sets.\\

{\bf The paper:}
In the first part of the present paper we study the counterpart of some aspects of Restle's paper in the fuzzy sets theory setting, when indicator functions of sets are replaced
by membership functions, whose definition we now recall (see for instance \cite{fuzz,zbMATH03226241}):

\begin{definition}
\label{seoul}
 A function from $X$ into $[0,1]$, i.e. belonging to  $[0,1]^X$, is called a membership function.
\end{definition}
We write a membership function $f$ as $\mu_{\widetilde{A}}$, where by definition, $\widetilde{A}$ denotes the fuzzy set defined by $f$.\smallskip

In machine learning, a recent research trend consists in replacing the real numbers by hypercomplex numbers; see for instance
\cite{aizenberg2013multi, kuroe} for complex numbers, \cite{MR4537587} for bicomplex numbers and \cite{kobayashi2012hyperbolic} for hyperbolic numbers.
In the second part of this paper, and inspired by the work \cite{MR3651492} where probabilities are allowed to take values in the set of hyperbolic numbers
(we will say for short, hyperbolic-valued), we initiate a
study of fuzzy set theory when the
membership function is hyperbolic-valued. Definitions are recalled in the sequel, but we already mention at this stage that hyperbolic numbers can be seen as the set of matrices of the form $\begin{pmatrix}a&b\\ b&a\end{pmatrix}$, where $a$ and $b$ run through the real numbers.\smallskip

We therefore address
two different audiences, the fuzzy set community and the hypercomplex analysis community, and will review materials from both fields to make the paper accessible to both groups.\smallskip

To pursue we recall that for two (not necessarily Hermitian) matrices $A$ and $B$ in $\mathbb C^{n\times n}$ one says that $A\le B$ if $B-A$ is a positive semi-definite matrix
(one also says non-negative), i.e. if $B-A$  is Hermitian (symmetric in the case of matrices with real entries) with non-negative eigenvalues.\smallskip

Two Hermitian matrices which do not commute cannot be simultaneously diagonalized and one cannot define in a natural way their maximum and minimum using the
natural order of matrices. On the other
hand, hyperbolic numbers are simultaneously diagonalizable and they form a lattice: we can define maximum and minimum (with respect to the above partial order)
of any pair of hyperbolic numbers in the set of hyperbolic numbers. As a consequence we can extend to the hyperbolic setting important
operations on fuzzy sets which involve maximum and minimum.\smallskip

As mentioned above in fuzzy set theory one replaces indicator functions of subsets of a given set $X$ by functions from $X$ into $[0,1]$.
We introduce a new operator on membership functions: given $f$ and $g$ two membership functions we associate the hyperbolic-valued function
\begin{equation}
\label{le-bossu}
  M_{f,g}(x)=\frac{1}{2}\begin{pmatrix}f(x)+g(x)&f(x)-g(x)\\
    f(x)-g(x)&f(x)+g(x)\end{pmatrix},\quad x\in X.
 \end{equation}

 Formula \eqref{le-bossu} defines a new operation on membership functions, and $M_{f,g}$ takes values in the counterpart of $[0,1]$ for hyperbolic numbers.
 The main properties of  this operation are obtained using the fact that the hyperbolic numbers form a lattice.\smallskip

 We note (see Section \ref{secfuzz} for definitions) that already in classical fuzzy set theory, fuzzy sets have been generalized to sets defined
 by two membership functions (intuitionistic fuzzy sets, also known as bipolar fuzzy sets, and soft fuzzy sets). The present extension
 is different from these approaches. \smallskip

 
    The hyperbolic  numbers form a commuting family of Hermitian matrices, and as such is simultaneously diagonalizable, as is also immediately seen from
    \eqref{cond1}. More generally recall that a commuting family of complex matrices is simultaneously triangularizable; see \cite{MR632835}.
    The present theory could be extended to families of commuting  symmetric matrices, or   diagonalizable familes of non-symmetric matrices.
\smallskip

The paper consists of ten sections besides the introduction. In Section 2 we discuss distances associated to positive definite kernels.
In Section 3 we discuss
betweenness of vectors in a vector space. Betweenness of sets is studied in Section 4. A few facts from fuzzy set theory are reviewed in Section 5.
Betweenness in the fuzzy setting is studied in Sections 6 and 7, using two different approaches: characterization in terms of intervals and in terms of strong $\alpha$-cuts. That the two defintions are equivalent is proved in Theorem \ref{bretagne}.
The definition and main properties of hyperbolic numbers are reviewed in Section 8, while
hyperbolic-valued membership functions are studied in Section 9 and their properties in Section 10. Betweenness in the setting of hyperbolic-valued
membership functions is considered in Section 11. \\

Finally, a word on notation: $a\wedge b$ and $a\vee b$ denote respectively the minimum and maximum of the real numbers $a$ and $b$, and more genrally the
corresponding operations in a lattice. The matrices
\[
  \begin{pmatrix}0&0\\0&0\end{pmatrix}\quad{\rm and}\quad\begin{pmatrix}1&0\\0&1\end{pmatrix}
\]
will be denoted sometimes by $0$ for the first and by $1$ or $I_2$ for the second.
\section{Positive definite kernel and associated metric}
\setcounter{equation}{0}{

We review some facts on positive definite functions relevant to the present work; for further references we suggest
\cite{MR2002b:47144,aron, saitoh}. For completeness we recall:

\begin{definition}
  \label{metric-space}
  Let $E$ be a set. The map $d$ from $E\times E$ into $[0,\infty)$ is called a metric (or a distance) if the following three conditions hold for all $x,y,z\in E$:
  \begin{eqnarray}    \label{sym5}
    d(x,y)&=& 0\quad\iff\quad x=y\\
    \label{sym6}
    d(x,y)&=&d(y,x)\\
    d(x,y)&\le&d(x,z)+d(y,z).
                \label{triangu}
\end{eqnarray}
\end{definition}

\eqref{triangu} is called the triangle inequality and the pair $(E,d)$ (or $E$ for short) is called a metric space.\\

Positive definite kernels (we will also say positive definite functions, although the latter terminology is usually used for a smaller class of kernels) whose definition
we now recall, play an important role in machine learning, in particular in the theory of support vector machines;
see \cite{pereverzyev2022introduction} for a recent account. Here they are of special interest because of the metric induced on the set where such a function is defined; 
see \cite{MR4302453} and see \cite{MR2002b:47144,aron, saitoh} for more information on positive definite kernels.

\begin{definition}
Let $E$ be a set and let $K(t,s)$ be defined on $E\times E$. It is called positive definite on $E$ if  for every choice of $N\in\mathbb N$ and $t_1,\ldots, t_N\in E$ the matrix $\left(K(t_\ell,t_j)\right)_{\ell,j=1}^N$ is positive semi-definite.
\end{definition}

The following classical theorem gives a characterization of positive definite functions; one direction is quite clear and in the other one can take $\mathcal H$ to be the reproducing kernel Hilbert space $\mathcal H(K)$ with reproducing kernel $K(t,s)$ since
\[
K(t,s)=\left\langle K(\cdot, s),K(\cdot, t)\right\rangle_{\mathcal H(K)}.
  \]

\begin{theorem}
  The function $K(t,s)$ is positive definite on $E$ if and only if it factors via a Hilbert space, i.e. if and only if there exists a Hilbert space
  $\mathcal H$ and a function $f_t$ from $E$ into $\mathcal H$ such that
  \begin{equation}
K(t,s)=\langle f_s,f_t\rangle_{\mathcal H},\quad t,s\in E.
\label{facto}
\end{equation}
  \end{theorem}

For the following proposition see for instance \cite{MR4302453}, where some explicit examples are also computed.
  
\begin{proposition}
Let $K(t,s)$ be positive definite on $E$ with factorization \eqref{facto}, and  assume that
\begin{equation}
  \label{yhn}
  t\not =s\quad\iff\quad f_t\not=f_s.
\end{equation}
Then,
\begin{equation}
  \label{voltaire}
d_K(t,s)=\sqrt{K(t,t)+K(t,s)-2{\rm Re}\, K(t,s)}
\end{equation}
defines a distance on $E$.
\end{proposition}

\begin{proof} 
  In view of \eqref{facto} we have
  \begin{equation}
    \label{45rt}
    d_K(t,s)=\|f_t-f_s\|_{\mathcal H}.
  \end{equation}
  and the three conditions for a metric follow, the first one using \eqref{yhn}. 
  \end{proof}

\begin{remark} For $f_t=K(\cdot, t)$ the condition \eqref{yhn} becomes
  \[
t\not=s\quad\iff\quad K(\cdot,s)\not\equiv K(\cdot,t).
    \]
    which is in fact a necessary and sufficient condition for \eqref{voltaire} to define a metric.
    \end{remark}

  As an example of metric $d_K$, let $\Omega$ be a finite set and let ${\rm Card}$ denote the counting measure. Then, the function
\[
  K(A,B)={\rm Card}\, (A\cap B)
\]
is positive definite on $\mathcal P(\Omega)$, and the associated metric is given by
\begin{equation}
  d(A,B)=\sqrt{{\rm Card}\, A+{\rm Card}\, B-{\rm Card}\, (A\cap B)}=\sqrt{{\rm Card}\, A\Delta B},\quad A,B\in\mathcal P(\Omega).
\end{equation}

In general the square of a metric is not a metric, but the squareroot of a metric is still a metric. In the present case, it so happens that the square of $d(A,B)$, namely
\begin{equation}
d^2(A,B)={\rm Card}\, A+{\rm Card}\, B-{\rm Card}\, (A\cap B)={\rm Card}\, A\Delta B,\quad A,B\in\mathcal P(\Omega),
\end{equation}
is still a metric (not induced by a positive definite kernel); see Proposition \ref{propoDm}. A weighted form of $d^2$ appear already in  \cite[p. 290-291]{hays1958approach} in the study of the difference
(called in \cite{hays1958approach} implicational difference) between traits in an individual. The distance $d^2(A,B)$ play a key role in the present work,
and an important difference between $d$ and $d^2$ will be shown in the paper.
\section{Betweenness of vectors}
\setcounter{equation}{0}
Let $\mathcal V$ be a real or complex vector space, and let $u,v\in\mathcal V$. Recall that the interval defined by $u$ and $v$ is the set of vectors of the form
\begin{equation}
  \label{gamma-t}
c(t)=u+t(v-u),\quad t\in[0,1]
  \end{equation}
  which reduces  to one point when $u=v$ (no order is assumed, and we can speak for instance of the interval $[0,1]$ as well as the interval $[1,0]$). It is natural to define:

  \begin{definition}
    \label{defdel}
The  vector $w$ is said to be between $u$ and $v$ if $w\in[u,v]$.
  \end{definition}

  \begin{remark}
    When extra structure is given on $\mathcal V$, or for algebraic structures different from a vector space structure, the above definition need not be possible, or even if possible, need not be the best one. For instance, in case of a lattice, a natural definition would be to replace \eqref{gamma-t} by
    \[
u\wedge v+t(u\vee v-u\wedge v), \quad t\in[0,1].
\]
When moreover a commutative product is available (as in the case of the hyperbolic numbers), one can replace $[0,1]$ by its counterpart with respect to the partial order; see Definition \ref{london}.
    \end{remark}
  
  Recall that a norm $\|\cdot\|$ on a vector space defines a metric via
  \[
d(u,v)=\|u-v\|.
    \]
  
    \begin{proposition}
      \label{equal-triangle}
Let $(\mathcal V,\|\cdot\|)$ be a normed space and let $w$ be between $u$ and $v$. Then equality holds in the triangle inequality for $d$, i.e.
\begin{equation}
  \label{d-eq}
      d(u,v)=d(u,w)+d(w,v).
      \end{equation}

\end{proposition}

  \begin{proof}
We write $w=c(t)$ where $t\in[0,1]$. We have:
    \[
      \begin{split}
        d(u,v)&=\|u-v\|\\
        d(u,w)&=\|u-(u-t(v-u))\|\\
        &=t\|u-v\|\\
        d(w,v)&=\|u+t(v-u)-v\|\\
        &=\|(1-t)(u-v)\|\\
        &=(1-t)\|u-v\|
      \end{split}
    \]
    and hence \eqref{d-eq} holds.
    \end{proof}

    The converse to the above claim is false in general, as can be seen by the example $\mathcal V=\mathbb R^2$ endowed with the norm $\|(x,y)\|_\infty
    =|x|\vee|y|.$ Take
    \[
      u=(0,0),\quad v=(1,1/4)\quad{\rm and}\quad w=(1/2,1/4).
    \]
    Then,
    \[
      \|u-v\|_\infty=\|u-w\|_\infty+\|w-v\|_\infty
    \]
    but $w\not\in[u,v]$. \smallskip

    The problem in the preceding example is that the norm is not strictly convex. We give the definition for complex vector spaces, but the same will hold for real vector spaces.

    \begin{definition}
      The norm $\|\cdot\|$ on the real or complex vector space $\mathcal V$  is called strictly convex if the following hold:
      \[
        \|u+v\|=\|u\|+\|v\|\quad and \quad u\not=0\Longrightarrow\quad \text{$v=cu$ for some $c\ge 0$.}
        \]
      \end{definition}

\begin{proposition}
  Assume the norm $\|\cdot\|$ strictly convex. Then
  \[
\|u-v\|=\|u-w\|+\|w-v\|\quad\iff\quad w\in[u,v]
    \]
  \end{proposition}

  \begin{proof}
    If $u=w$ the result is trivial. Assume therefore $u\not=w$.
    Since $u-v=(u-w)+(w-v)$ it follows from the definition that $w-v=c(u-w)$ for some $c\ge 0$. Hence
    \[
w=\frac{c}{1+c}u+\frac{1}{1+c}v=u+t(v-u)
\]
with\footnote{$t=0$ corresponds to $c\longrightarrow\infty$} $t=\frac{1}{1+c}\in[0,1]$.
    \end{proof}

Examples of norms on the  dimensional vector space $\mathbb C^N$ are given by (with $z=(z_1,\ldots, z_N)$ and similarly for $w$)
\[
  \|z\|_r=\begin{cases}\,\left(\sum_{n=1}^Nm_n|z_n|^r\right)^{1/r},\quad r\in[1,\infty)\\
   \, \vee_{n=1}^Nm_n|z_n|,\,\quad\hspace{1.4cm} r=\infty
    \end{cases}
  \]
  where $m_1,\ldots, m_N$ are strictly positive. They are strictly convex for $1<r<\infty$ but not for $r\in\left\{1,\infty\right\}$, and correspond
to the distances

\[
  D_r(z,w)=\begin{cases}\,\left(\sum_{n=1}^Nm_n|z_n-w_n|^r\right)^{1/r},\quad r\in[1,\infty)\\
   \, \vee_{n=1}^Nm_n|z_n-w_n|,\,\quad\hspace{1.4cm} r=\infty.
    \end{cases}
  \]

These norms fall into a larger family of norms used  \cite{zwick1987measures} for membership functions, and which may be defined as follows. We will assume that
$(X,\mathcal A,\sigma)$ is a measured space, with sigma-algebra $\mathcal A$ and positive measure $\sigma$.
The measure $\sigma$ has the following properties (which allows to define these norms on membership functions, since the latter take values in $[0,1]$)

\begin{definition}
$\sigma$ will be a positive measure such that $\int_Xd\sigma(x)<\infty$ and with the condition:
\begin{equation}
  \label{one-to-one}
\int_X|f(x)|d\sigma(x)=0\quad \Longrightarrow \quad f=0,\,\, a.e.
\end{equation}
\label{defsigmaa}
\end{definition}

We will say that two sets in $\mathcal A$ are equivalent (notation: $A\sim B$) if $\sigma(A\Delta B)=0$. We have an equivalent relation since:\\
$(1)$ It is reflexive since $A\Delta A=\emptyset$ and so $\sigma(A\Delta A)=0$.\\
$(2)$ It is symmetric since $A\Delta B=B\Delta A$.\\
$(3)$ It is transitive. For $A,B,C\in\mathcal N$ assume $A\sim B$ and $B\sim C$. Then $A\sim C$ since
\[
  A\Delta C=(A\Delta B)\Delta (B\Delta C)
\]
and
\[
  \sigma(A\Delta C)=\sigma((A\Delta B)\Delta (B\Delta C))\le \sigma(A\Delta B)+\sigma(B\Delta C)=0.\\
\]
\begin{definition}
  We denote by $\mathcal N_0$ the  elements of $\mathcal A$ equivalent to $\emptyset$ and by $\mathcal A_0=\mathcal A/\mathcal N$ the space of equivalent classes.
  \label{equiv12}
\end{definition}

We set, for $f$ measurable and bounded in modulus,
\[
\|f\|_r=\left(\int_X|f(x)|^rd\sigma(x)\right)^{1/r},\quad r\in[1,\infty)
\]
and
  \[
 \|f\|_\infty=\esssup_{x\in X}|f(x)|,
\]
corresponding to distances $d_r$ and $d_\infty$.\smallskip

The following result is of limited interest since the number $t$ in \eqref{t-ind-x} does not depend on $x$, but stresses the difference with the results presented in
Sections \ref{secfuzz} and \ref{sec-alpha}.

\begin{proposition} Given $(X,\mathcal A,\sigma)$ a measured space,
  assume that $\mu_{\widetilde{A}},\mu_{\widetilde{B}}$ and $\mu_{\widetilde{C}}$ are measurable membership functions, and that
  $\mu_{\widetilde{C}}$ is between $\mu_{\widetilde{A}}$ and $\mu_{\widetilde{C}}$ in the sense of Definition \ref{defdel} meaning that
there exists $t\in[0,1]$, {\sl independent of $x$}, such that
          \begin{equation}
            \label{t-ind-x}
            \mu_{\widetilde{C}}(x)= t\mu_{\widetilde{A}}(x)+(1-t)\mu_{\widetilde{B}},\quad x\in X.
                  \end{equation}
  Then,
  \begin{equation}
d_r(\mu_{\widetilde{A}},\mu_{\widetilde{B}})=d_r(\mu_{\widetilde{A}},\mu_{\widetilde{C}})+d_r(\mu_{\widetilde{C}},\mu_{\widetilde{B}}),\quad \forall r\in[1,\infty]
\end{equation}
The converse statement is true if $r\not\in\left\{1,\infty\right\}$.
  \end{proposition}

\begin{proof}
  The direct claim follows from Proposition \ref{equal-triangle}.
    We now turn to the converse statement. Since $r\in(1,\infty)$ the norm $d_r$ is strictly convex. Thus equality in the triangle inequality means that $\mu_{\widetilde{C}}$ is in
    the interval  defined by $\mu_{\widetilde{A}}$  and $\mu_{\widetilde{B}}$.
\end{proof}

\section{Betweenness of sets}
\setcounter{equation}{0}
In preparation for the following sections we rewrite in a slightly different form some  results from \cite{restle1959metric}.
We first recall a definition.

\begin{definition}
Let $X$ be some non-empty set. 
A set $A\in\mathcal P(X)$ is uniquely determined by its indicator function $1_A$ defined by
\begin{equation}
  \label{indic}
  1_A(x)=\begin{cases}\, 1\,\,{\rm if}\,\, x\in A\\
  \,0\,\,{\rm if}\,\,x\not\in A.
\end{cases}
\end{equation}
\end{definition}

There is therefore in classical set theory a one-to-one correspondence between elements of $\mathcal P(X)$ and
the set $\left\{0,1\right\}^X$ of functions from $X$ into $\left\{0,1\right\}$.
As is well known, the indicator functions of the union, intersection and symmetric difference of two sets $A$ and $B$ and of the complement of a
set $A$ are given by
\begin{eqnarray}
  \label{llg1}
  1_{A\cup B}&=&1_A+1_B-1_A1_B\\
  &=&1_A\vee1_B\\
      1_{A\cap B}&=&1_A\wedge1_B\\
  &=&1_A1_B\\
  1_{A\Delta B}&=&1_A+1_B-2\cdot 1_A1_B\\
  \label{llg6}
             &=&1_A\vee1_B-1_A\wedge 1_B\\
  &=&(1_A-1_B)^2\\
  1_{X\setminus A}&=&1-1_A
                  \label{llg7}
\end{eqnarray}
where we have denoted by $\Delta$ the symmetric difference and by $X\setminus A$ the complement of the set $A$.

\begin{lemma} Let $X$ be a set and let $A,B,C\in\mathcal P(X)$.
  Then, $C$ is between $A$ and $B$ in the sense of equation \eqref{ineq456} if and only if
  \begin{equation}
    \label{paris1}
    1_C=1_{A\cap B}+1_Z
  \end{equation}
  where $Z\in\mathcal P(X)$ is such that
  \begin{equation}
    \label{paris2}
    Z\subset A\Delta B.
  \end{equation}
  \label{lemma4-2}
\end{lemma}

\begin{proof}
  We rewrite \eqref{crispbetween} as \eqref{ineq456}, i.e.
  \[
     1_A(x)\wedge 1_B(x)\le 1_C(x)\le 1_A(x)\vee1_B(x),\quad \forall x\in\mathbb R.
\]
Thus there exists $t(x)\in[0,1]$ such that
\[
1_C(x)=1_A(x)\wedge1_B(x)+t(x)\underbrace{\left( 1_A(x)\vee 1_B(x)-1_A(x)\wedge1_B(x)\right)}_{1_{A\Delta B}(x)}.
\]
From the above, $t(x)$ can be chosen to belong to $\left\{0,1\right\}$.
Define a set $Z$ via
\[
1_Z(x)=t(x)\underbrace{\left(1_A(x)\vee1_B(x)-1_A(x)\wedge1_B(x)\right)}_{1_{A\Delta B}(x)}.
\]
Then $Z\subset A\Delta B$. The converse is clear.
\end{proof}
Given two sets $A$ and $B$ in $\mathcal P(\Omega)$ we note that the interval $[1_A,1_B]$ is not made of indicator functions in general,
but consists of the functions of the form
\begin{equation}
  \label{member1}
f_t(x)=t1_A(x)+(1-t)1_B(x)
\end{equation}
when $t$ varies in $[0,1]$. A more interesting case is when $t$ is allowed to vary with $x$:
\begin{equation}
  \label{member1}
f_t(x)=t(x)1_A(x)+(1-t(x))1_B(x),
\end{equation}
where now $t$ is a function from $X$ into $[0,1]$.
These functions are examples of membership functions, which are the main tool in fuzzy set theory. The interpretation of \eqref{member1} in terms of norms uses
the $\alpha$-cuts. See Section \ref{sec-alpha} below.\smallskip


For the following result, see also Restle \cite{restle1959metric}, where the importance of the equality case in the triangle inequality is stressed out.

\begin{proposition}
  Let $(X,\mathcal A,\sigma)$ be as above with $\sigma$ a measure on $X$ satisfying the hypothesis of Definition \ref{defsigmaa}, and let $\mathcal A_0$ be as in Definition \ref{equiv12}. Then,
  \[
D_\sigma(A_0,B_0)=\int_X(1_A(x)-1_B(x))^2d\sigma(x)=\sigma(A\Delta B),\quad A,B\in\mathcal A_0
\]
(where $A\in\mathcal A$ belongs to the equivalence class of $A_0$ and $B\in\mathcal A$ belongs to the equivalence class of $B_0$)
is a metric on $\mathcal A_0$. If $X$ has finite cardinal one can take $\mathcal A_0=\mathcal P(X)$.
\label{propoDm}
\end{proposition}

\begin{proof}
  The various definitions do not depend on the chosen representative in a given equivalence class.
Assume that $D_\sigma(A_0,B_0)=0$. By \eqref{one-to-one} we have $\sigma(A\Delta B)=0$ and hence $A_0=B_0$.
  It is clear that $D_\sigma(A_0,B_0)=D_\sigma(B_0,A_0)$. We now check the triangle inequality and first note that  (with $C_0\in\mathcal A_0$ and $C$ in the equivalence class $C_0$)
    \begin{equation}
     (1_A-1_C)^2(x)+(1_C- 1_B)^2(x)-(1_A-1_ B)^2(x)=2(1_A-1_C)(x)(1_B- 1_C)(x)
     \label{equa-q}
   \end{equation}
   for $x\in X$.
   Hence, we have
   \begin{equation}
     \label{678p}
     \begin{split}
       D_\sigma(A_0,C_0)+D_\sigma(B_0,C_0)-D_\sigma(A_0,B_0)&=\\
       &\hspace{-3cm}=\int_X\left\{(1_A-1_C)^2(x)+(1_C- 1_A)^2(x)-(1_A-1_ B)^2(x)\right\}d\sigma(x)\\
       &\hspace{-3cm}=2\int_X(1_A-1_C)(x)(1_B- 1_C)(x)d\sigma(x).
     \end{split}
   \end{equation}

But
\begin{equation}
  \label{567890}
        (1_A-1_C)(x)(1_B- 1_C)(x=\begin{cases}\,\,\,\,\,\quad\hspace{0.6cm}(1-1_C(x))^2=1-1_C(x),\quad x\in A\cap B\\
          \,\, -1_C(x)(1-1_C(x))=0,\quad \quad \hspace{1.14cm}x\in  B\setminus A\\
          \,\,\,\,\,\,\,(1 -1_C(x))1_C(x)=0,\quad \quad \hspace{1.1cm}x\in  A\setminus B\\
          \qquad\hspace{1.6cm}1_C(x)^2=1_C(x),\quad\hspace{0.64cm} x\in X\setminus (A\cup B)
          \end{cases}
        \end{equation}
       and so
    \[
        \int_X(1_A-1_C)(x)(1_B- 1_C)(x)d\sigma(x)=\int_{A\cap B}(1-1_C(x))d\sigma(x)+\int_{X\setminus (A\cup B)}1_C(x)d\sigma(x)
     \]
Hence \eqref{678p} becomes
\[
  D_\sigma(A_0,C_0)+D_\sigma(B_0,C_0)-D_\sigma(A_0,B_0)=\int_{A\cap B}(1-1_C(x))d\sigma(x)+\int_{X\setminus (A\cup B)}1_C(x)d\sigma(x),
\]
which is non-negative, 
and hence the triangle inequality holds for $D_\sigma$.
\end{proof}

\begin{proposition}
  \label{ineq-sigma-eta}
  In the notation of the previous proposition,  $C_0$ is between $A_0$ and $B_0$ for the metric $D_\sigma(A_0,B_0)$ if and only if the triangle inequality is an equality:
  \begin{equation}
    \label{3456}
D_\sigma(A_0,B_0)=D_\sigma(A_0,C_0)+D_\sigma(C_0,B_0)
    \end{equation}
  \end{proposition}

  \begin{proof}
 By the triangle inequality for $D_\sigma$, and using \eqref{equa-q} we have:
     \[
      \begin{split}
        0&\le \int_X\left\{1_{A\Delta C}(x)+1_{C\Delta B}(x)-1_{A\Delta B}(x)\right\}d\sigma(x)\\
        &=\int_X\left\{(1_A-1_C)^2(x)-(1_C- 1_B)^2(x)-(1_A-1_ B)^2(x)\right\}d\sigma(x)\\
        &=2\int_X(1_A-1_C)(x)(1_B- 1_C)(x)d\sigma(x).
              \end{split}
      \]
      
 Hence,
        \[
          \begin{split}
            0&\le \int_X\left\{1_{A\Delta C}(x)+1_{C\Delta B}(x)-_{A\Delta B}(x)\right\}d\sigma(x)\\
            &=2\left\{      \int_{A\cap B}(1-1_C(x))d\sigma(x)+\int_{X\setminus (A\cup B)}1_C(x)\right\}d\sigma(x)
          \end{split}
          \]
    
    Hence, \eqref{3456} holds if and only if
    \[
      \int_{A\cap B}(1-1_C(x))d\sigma(x)=\int_{X\setminus (A\cup B)}1_C(x)d\sigma(x)=0
    \]
    that is, if and only if
    \[
      1_C(x)=1,\quad x\in A\cap B\quad{\rm and}\quad 1_C(x)=0, \quad x\not\in A\cup B,
    \]
    i.e. if and only if \eqref{crispbetween} is in force.
    \end{proof}

    $\sqrt{D_\sigma}$ is also  a metric on $\mathcal A$ (maybe more natural a priori since it arises from a positive definite kernel), but we have:
\begin{proposition}
 Let $A_0,B_0,C_0\in\mathcal A_0$. It holds that
  \begin{equation}
  \label{sqrt-eqqua}
\sqrt{D_\sigma(A_0,B_0)}=\sqrt{D_\sigma(A_0,C_0)}+\sqrt{D_\sigma(C_0,B_0)}
\end{equation}
if and only if $C_0=A_0$ or $C_0=B_0$.
\end{proposition}

\begin{proof}
  Assume that \eqref{sqrt-eqqua} is in force. Taking square and taking into account \eqref{equa-q} we obtain, with $A,B,C$ being in the equivalence classes of $A_0,B_0$ and $C_0$ respectively
  \begin{equation}
    \label{c-s-123}
-2\int_X(1_A-1_C)(x)(1_B- 1_C)(x)d\sigma(x)=2\sqrt{D_\sigma(A_0,C_0)}\sqrt{D_\sigma(C_0,B_0)}.
\end{equation}
We are thus in the equality case in the Cauchy-Schwarz inequality. If $C=A$ there is nothing to prove.  Assuming $C_0\not=A_0$,  there exists $u\in\mathbb C$ such that
\begin{equation}
  \label{A-B-C}
(1_B-1_C)=u(1_C-1_A),\quad \sigma\,\, a.e.
\end{equation}
Plugging this  into \eqref{c-s-123} we obtain
\[
u\int_X(1_C(x)-1_A(x))^2d\sigma(x)=|u|\cdot\int_X(1_C(x)-1_A(x))^2d\sigma(x).
  \]
  Since $A\not=C$ it follows that $u\ge 0$.\smallskip

  If $u=0$ in \eqref{A-B-C} we have $B=C$ and hence $B_0=C_0$. We now show by contradiction that we cannot have $u\not=0$ since $C\not=A$. Assume thus $u\not=0$ (and so $u>0$) and first suppose that there is $x\in C\setminus A$. Then, \eqref{A-B-C} becomes
  \[
    (1_B-1)=u.\]
  The left handside of this equality is less or equal to $0$ while the right handside is strictly positive, which is impossible. Assume now that there is $x\in A\setminus C$ Then \eqref{A-B-C} becomes
  \[
1_B=-u
    \]
  which is impossible for the same reason as above.
\end{proof}

  \section{Fuzzy set theory}
\setcounter{equation}{0}
\label{secfuzz}
In a way similar to information theory, which originates in 1948 with Shannon's paper \cite{MR10:133e}, one can pinpoint the origin of
fuzzy set theory and logic with the papers of Zadeh \cite{zbMATH03226241}, but it is good to mention the earlier works on multi-valued logic of
Lukasiewicz \cite{zbMATH03335866}.
For the convenience of the reader we review some definitions from fuzzy set theory, and send the reader to the books
\cite{chen-pham,MR1356718,fuzz,zbMATH01202379,MR1174743} for further information.

The set of indicator functions  is $\left\{0,1\right\}^X$, and is therefore included in the set of membership functions (see Definition \ref{seoul} for the latter).
Let $N\in\mathbb N$. We note that to any function from $[0,1]^N$ into $[0,1]$ one can define a map which to $N$ membership functions
associates a new membership function.\smallskip

Each of the functions \eqref{llg1}-\eqref{llg6} (and to a certain extent also \eqref{llg7}) have numerous possible extensions in the setting
of membership functions. This degree of freedom is one of the main strengths of fuzzy set theory. As a first example, consider the intersection, with indicator function $1_A1_B$.  When $1_A$ and $1_B$ are replaced by membership functions $\mu_{\widetilde{A}}$ and
$\mu_{\widetilde{B}}$, two different extensions of intersection of classical sets
will be given by $\mu_{\widetilde{A}}\vee\mu_{\widetilde{B}}$ and $\mu_{\widetilde{A}}\mu_{\widetilde{B}}$.\smallskip

The maximum and minimum of two membership functions $\mu_{\widetilde{A}}$ and $\mu_{\widetilde{B}}$ are also membership functions, corresponding respectively to the union
$\widetilde{A}\cup\widetilde{B}$  and intersection $\widetilde{A}\cap\widetilde{B}$ of the fuzzy sets $\widetilde{A}$  and $\widetilde{B}$  (see \cite[\S 3.1 p. 30]{fuzz}).\smallskip

The product $\mu_{\widetilde{A}}\mu_{\widetilde{B}}$ is also a membership functions, corresponding to a fuzzy set called algebraic product of the fuzzy sets $\widetilde{A}$ and $\widetilde{B}$; see \cite[\S 3.3 p. 33]{fuzz}).\smallskip

For $\mu_{\widetilde{A}}$ a membership function,  $1-\mu_{\widetilde{A}}$ is still a membership function, corresponding to a fuzzy set called the fuzzy complement of $\widetilde{A}$, and denote by $cl(\widetilde{A})$.\smallskip

More generally, one can take functions with values in a lattice; this was already done by Zadeh's student Goguen, see \cite{MR224391}, and
later also developped by Atanassov in his theory of intuitionistic sets; see  \cite{MR852871}. 
An intuitionistic fuzzy set defined on a set $X$ is defined by two functions from $X$ into $[0,1]$,
respectively called membership function and non-membership function. In the second part of this paper (Sections \ref{section2}-\ref{sec-end}) we will consider
the lattice of hyperbolic numbers.

\section{Betweenness in the fuzzy case}
\setcounter{equation}{0}

\begin{definition}
  Let $\mu_{\widetilde{A}},\mu_{\widetilde{B}}$ and $\mu_{\widetilde{C}}$ be membership functions.  We say that $\mu_{\widetilde{C}}$ is pointwise  between
  $\mu_{\widetilde{A}}$ and $\mu_{\widetilde{B}}$ if
  \begin{equation}
    \mu_{\widetilde{A}}(x)\wedge\mu_{\widetilde{B}}(x)
    \le\mu_{\widetilde{C}}(x)\le     \mu_{\widetilde{A}}(x)\vee \mu_{\widetilde{B}}(x),\quad \forall x\in X.
    \label{maxminmax}
    \end{equation}
  \end{definition}

  In other words, for every $x\in X$, $\mu_{\widetilde{C}}(x)$ belongs to the interval determined by $\mu_{\widetilde{A}}(x)$ and $\mu_{\widetilde{B}}(x)$.\\
  
As it should be this definition reduces to \eqref{crispbetween} in the crisp case since then we have
\[
1_{A\cap B}(x)=\mu_{\widetilde{A}}(x)\wedge\mu_{\widetilde{B}}(x)\quad{\rm and}\quad    1_{A\cup B}(x)=\mu_{\widetilde{A}}(x)\vee\mu_{\widetilde{B}}(x),\quad \forall x\in\mathbb R.
\]
The counterpart of Lemma \ref{lemma4-2} is as follows:
  \begin{proposition}
    Let $\mu_{\widetilde{A}},\mu_{\widetilde{B}}$ and $\mu_{\widetilde{C}}$ be membership functions. Then,  $\mu_{\widetilde{C}}$ is pointwise between $\mu_{\widetilde{A}}$ and $\mu_{\widetilde{B}}$ if and only if
  \begin{equation}
    \label{nocrsip1}
    \mu_{\widetilde{C}}(x)=\mu_{\widetilde{A}}(x)\wedge\mu_{\widetilde{B}}(x)+\mu_{\widetilde{Z}}(x)
  \end{equation}
  where $\mu_{\widetilde{Z}}$ is a membership function satisfying
  \begin{equation}
    \label{nocrsip2}
\mu_{\widetilde{Z}}(x)\le \mu_{\widetilde{A}}(x)\vee\mu_{\widetilde{B}}(x) -   \mu_{\widetilde{A}}(x)\wedge\mu_{\widetilde{B}}(x).\\
  \end{equation}
  \label{prop6-2}
  \end{proposition}

  \begin{proof}
  Assume first that \eqref{nocrsip1} and \eqref{nocrsip2} are in force. \eqref{nocrsip1} implies that
  \[
    \mu_{\widetilde{C}}(x)\ge \mu_{\widetilde{A}}(x)\wedge\mu_{\widetilde{B}}(x),
  \]
  and  \eqref{nocrsip1} and \eqref{nocrsip2} together lead to
  \[
    \begin{split}
      \mu_{\widetilde{C}}(x)&\le \mu_{\widetilde{A}}(x)\wedge\mu_{\widetilde{B}}(x)
      +\mu_{\widetilde{A}}(x)\vee\mu_{\widetilde{B}}(x) -   \mu_{\widetilde{A}}(x)\wedge\mu_{\widetilde{B}}(x)\\
      &=
      \mu_{\widetilde{A}}(x)\vee\mu_{\widetilde{B}}(x).
      \end{split}
    \]
    Conversely, assume that \eqref{maxminmax} holds. The formula
    \[
\mu_{\widetilde{Z}}(x)=\mu_{\widetilde{C}}(x)-\mu_{\widetilde{A}}(x)\wedge\mu_{\widetilde{B}}(x)
\]
defines a membership function which answers the question.
\end{proof}

\begin{remark}
  In the crisp case, \eqref{nocrsip1}-\eqref{nocrsip2} reduce to \eqref{paris1}-\eqref{paris2}.
  \end{remark}

  \section{Betweenness in the fuzzy case with $\alpha$-cuts}
  \setcounter{equation}{0}
\label{sec-alpha}
  Metrics between membership functions using $\alpha$-cuts have been defined in \cite{ralescu1984probability}.
  For $\alpha\in [0,1]$ we consider the strong $\alpha$-cuts $A_\alpha^\prime$ associated to the membership function $\mu_{\widetilde{A}}$, defined by
  \begin{equation}
    A_\alpha^\prime=     \mu_{\widetilde{A}}^{-1}(\alpha,1],\quad \alpha\in[0,1].
  \end{equation}
  Note than one also defines $\alpha$-cuts
    \begin{equation}
    A_\alpha=     \mu_{\widetilde{A}}^{-1}[\alpha,1],\quad \alpha\in[0,1].
  \end{equation}
  See e.g. \cite[p. 14]{MR1174743}.
The arguments in this section will not hold with the latter definition; strict inequalities are needed.
  
  \begin{definition}
    \label{orleans}
    Let $\mu_{\widetilde{A}},\mu_{\widetilde{B}}$ and $\nu_{\widetilde{C}}$ be membership functions from $X$ to $[0,1]$. We say that $\mu_{\widetilde{C}}$ is $\alpha$-between $\mu_{\widetilde{A}}$ and $\mu_{\widetilde{B}}$ if
$\mu_{\widetilde{C}}^{-1}(\alpha,1]$ is between $\mu_{\widetilde{A}}^{-1}(\alpha,1]$  and $\mu_{\widetilde{B}}^{-1}(\alpha,1]$  for every $\alpha\in[0,1]$.
    \end{definition}

    \begin{theorem}
      \label{bretagne}
      $\mu_{\widetilde{C}}$ is $\alpha$-between $\mu_{\widetilde{A}}$ and $\mu_{\widetilde{B}}$
      if and only if $\mu_{\widetilde{C}}$ is pointwise between $\mu_{\widetilde{A}}$ and $\mu_{\widetilde{B}}$.
\end{theorem}

    \begin{proof}
      We first assume that $\mu_{\widetilde{C}}$ is $\alpha$-between $\mu_{\widetilde{A}}$ and $\mu_{\widetilde{B}}$.
      Let $x\in X$ and let $\mu_{\widetilde{A}}(x),\mu_{\widetilde{B}}(x)$ and $\mu_{\widetilde{C}}(x)$ be the corresponding values of the
      membership functions. Since $\mu_{\widetilde{A}}
      $ and $\mu_{\widetilde{B}}$ play a symmetric role we can assume without loss of generality that
      \begin{equation}
        \label{1qaz}
        \mu_{\widetilde{A}}(x)\le\mu_{\widetilde{B}}(x).
      \end{equation}
      We want to show that \eqref{maxminmax} holds for all $x\in X$,  i.e.  taking into account \eqref{1qaz}, that
      \begin{equation}
        \label{max2}
        \mu_{\widetilde{A}}(x)\le\mu_{\widetilde{C}}(x)\le \mu_{\widetilde{B}}(x).
       \end{equation}
Equivalently we have to show that the following cannot hold:
      \begin{equation}
        \label{impo1}
        \mu_{\widetilde{C}}(x)<\mu_{\widetilde{A}}(x)
      \end{equation}
      or
      \begin{equation}
         \mu_{\widetilde{B}}(x)<\mu_{\widetilde{C}}(x).
\label{impo2}
      \end{equation}

      Assume first by contradiction that \eqref{impo1} holds. So $x\in
      \mu_{\widetilde{A}}^{-1}(
      \mu_{\widetilde{C}}(x),1] $. By \eqref{1qaz} we also have $x\in\mu_{\widetilde{B}}^{-1}(\mu_{\widetilde{C}}(x),1]$, and so
    \[
\mu_{\widetilde{A}}^{-1}(\mu_{\widetilde{C}}(x),1]\cap \mu_{\widetilde{B}}^{-1}(\mu_{\widetilde{C}}(x),1]\not=\emptyset
\]
since $x$ belongs to the intersection.
But
\[
\mu_{\widetilde{C}}^{-1}(\mu_{\widetilde{C}}(x),1]=\left\{y\in X\,\, :\,\, \mu_{\widetilde{C}}(x)<\mu_{\widetilde{C}}(y)\le 1\right\}
  \]
  and so $x\not \in\mu_{\widetilde{C}}^{-1}(\mu_{\widetilde{C}}(x),1]$, leading to a contradiction since, with $\alpha=\mu_{\widetilde{C}}(x)$ the hypothesis of $\alpha$-betweenness of $\mu_{\widetilde{C}}$ between $\mu_{\widetilde{A}}$ and $\mu_{\widetilde{B}}$
  gives
  \[
    \underbrace{\mu_{\widetilde{A}}^{-1}(\mu_{\widetilde{C}}(x),1]\cap \mu_{\widetilde{B}}^{-1}(\mu_{\widetilde{C}}(x),1]}_{{\rm contains}\,\, x}
  \subset\underbrace{\mu_{\widetilde{C}}^{-1}(\mu_{\widetilde{C}}(x),1].}_{{\rm does \,\, not\,\, contain}\,\, x}
    \]

      Assume now by contradiction that \eqref{impo2} holds. Then $x\in \mu_{\widetilde{C}}^{-1}(\mu_{\widetilde{B}}(x),1]$ but
  $x\not\in\mu_{\widetilde{B}}^{-1}(\mu_{\widetilde{B}}(x),1]$. The $\alpha$-betweenness with $\alpha=\mu_{\widetilde{B}(x)}$ gives the
 the decomposition
  \[
    \mu_{\widetilde{C}}^{-1}(\mu_{\widetilde{B}}(x),1]=
    \underbrace{\left(\mu_{\widetilde{A}}^{-1}(\mu_{\widetilde{B}}(x),1]\cap \mu_{\widetilde{B}}^{-1}(\mu_{\widetilde{B}}(x),1]\right)}_{\mbox{
            $\emptyset$ since $x\not \in  \mu_{\widetilde{B}}^{-1}(\mu_{\widetilde{B}}(x),1]$} }    \cup\, Z_x
        \]
        with
        \[
        Z_x\subset \left(\mu_{\widetilde{A}}^{-1}(\mu_{\widetilde{B}}(x),1]\setminus\mu_{\widetilde{B}}^{-1}(\mu_{\widetilde{B}}(x),1]\right)\cup
        \underbrace{ \left(\mu_{\widetilde{B}}^{-1}(\mu_{\widetilde{B}}(x),1]\setminus\mu_{\widetilde{A}}^{-1}(\mu_{\widetilde{B}}(x),1]\right).}_{\mbox{
            $\emptyset$  since $x\not \in  \mu_{\widetilde{B}}^{-1}(\mu_{\widetilde{B}}(x),1]$}}
          \]
See Lemma \ref{lemma4-2} and equation \eqref{paris2}. Since $x\not \in  \mu_{\widetilde{B}}^{-1}(\mu_{\widetilde{B}}(x),1]$ we have that $x\in Z_x$ and in particular  $x\in \mu_{\widetilde{A}}^{-1}(\mu_{\widetilde{B}}(x),1]$, so that
                  \begin{equation}
                    \mu_{\widetilde{A}}(x)>\mu_{\widetilde{B}}(x),
                  \end{equation}
                  contradicting \eqref{1qaz}.\\
                
                  Conversely, assume that $\mu_{\widetilde{C}}$ is pointwise between $\mu_{\widetilde{A}}$ and $\mu_{\widetilde{B}}$.
                  Thus, for every $x\in X$,
  \begin{equation}
    \mu_{\widetilde{C}}(x)\in[\mu_{\widetilde{A}}(x)\wedge\mu_{\widetilde{B}}(x),\mu_{\widetilde{A}}(x)\vee\mu_{\widetilde{B}}(x)]
  \label{inb}
\end{equation}
is in force. We want to show that, for every $\alpha\in[0,1]$
\begin{equation}
  \label{7890}
  \mu_{\widetilde{A}}^{-1}(\alpha,1]\cap\mu_{\widetilde{B}}^{-1}(\alpha,1]
  \subset \mu_{\widetilde{C}}^{-1}(\alpha,1]\subset
  \mu_{\widetilde{A}}^{-1}(\alpha,1]\cup\mu_{\widetilde{B}}^{-1}(\alpha,1].
\end{equation}
We divide this part of the proof in a number of steps.\\

      STEP 1: {\sl If there is no $x$ such that $\mu_{\widetilde{C}}(x)>\alpha$, both inclusions in \eqref{7890}
      are satisfied.}\smallskip

      The second inclusion in \eqref{7890} is now trivial. We show that the first one holds (and reduces to $\emptyset=\emptyset$).
      By hypothesis,
      \begin{equation}
        \label{567yhn}
        \mu_{\widetilde{C}}(x)\le \alpha
\end{equation}
      for all $x\in X$. Assume by contradiction
      that there is $y\in \mu_{\widetilde{A}}^{-1}(\alpha,1]\cap\mu_{\widetilde{B}}^{-1}(\alpha,1]$. Then
          \[
          \alpha<\mu_{\widetilde{A}}(y)\le 1\quad{\rm and}\quad           \alpha<\mu_{\widetilde{B}}(y)\le 1.
          \]
          In particular
          \begin{equation}
            \label{1234}
          \alpha<\mu_{\widetilde{A}}(y)\wedge\mu_{\widetilde{B}}(y).
          \end{equation}
          By the hypothesis on pointwise betweenness
          \begin{equation}
            \label{56789}
\mu_{\widetilde{A}}(y)\wedge\mu_{\widetilde{B}}(y)\le\mu_{\widetilde{C}}(y)
          \end{equation}
          Equations \eqref{567yhn}, \eqref{1234} and \eqref{56789} lead to
          \[
          \alpha<\mu_{\widetilde{A}}(y)\wedge\mu_{\widetilde{B}}(y)\le\mu_{\widetilde{C}}(y)\le\alpha,
          \]
          which cannot be.\\

          STEP 2: {\sl The first inclusion in \eqref{7890} holds.}\smallskip

                    If $\mu_{\widetilde{A}}^{-1}(\alpha,1]\cap\mu_{\widetilde{B}}^{-1}(\alpha,1]=\emptyset$ the first inclusion is trivially met.
Assume now that there is  $x\in\mu_{\widetilde{A}}^{-1}(\alpha,1]\cap\mu_{\widetilde{B}}^{-1}(\alpha,1].$
    Then, $x$ is such that $\mu_{\widetilde{A}}(x)\in(\alpha,1]$ and $\mu_{\widetilde{B}}(x)\in(\alpha,1]$. Thus
        \[
\alpha<\mu_{\widetilde{A}}(x)\le 1\quad{\rm and}\quad \alpha<\mu_{\widetilde{B}}(x)\le 1.
        \]
        From $\mu_{\widetilde{A}}(x)\wedge\mu_{\widetilde{B}}(x)\le \mu_{\widetilde{C}}(x)$ we have that
        $x\in\mu_{\widetilde{C}}^{-1}(\alpha,1]$, and hence
  \[
\mu_{\widetilde{A}}^{-1}(\alpha,1]\cap\mu_{\widetilde{B}}^{-1}(\alpha,1]\subset\mu_{\widetilde{C}}^{-1}(\alpha,1].
      \]

STEP 3: {\sl The second inclusion in \eqref{7890} holds.}\smallskip

          By Step 1 we may assume that $\mu_{\widetilde{C}}^{-1}(\alpha,1]\not=\emptyset$. Let thus $x$ be such that
$\mu_{\widetilde{C}}(x)>\alpha$.
Since $\mu_{\widetilde{A}}(x)\vee\mu_{\widetilde{B}}(x)\ge \mu_{\widetilde{C}}(x)$. we have
\[
\mu_{\widetilde{C}}^{-1}(\alpha,1]\subset
\mu_{\widetilde{A}}^{-1}(\alpha,1]\cup\mu_{\widetilde{B}}^{-1}(\alpha,1].
\]
Hence $x\in\mu_{\widetilde{A}}^{-1}(\alpha,1]\cup\mu_{\widetilde{B}}^{-1}(\alpha,1]$ and thus $\mu_{\widetilde{C}}(\alpha,1]$ is
between $\mu_{\widetilde{A}}^{-1}(\alpha,1]$ and $\mu_{\widetilde{B}}^{-1}(\alpha,1]$.
\end{proof}

Let now $\sigma$ be a measure on $X$ satisfying the properties of Definition \ref{defsigmaa} and let $\eta$ be a strictly positive measure on $[0,1]$
We define (assuming the integrals well defined)

\begin{equation}
  \label{DMM}
  D(\mu_{\widetilde{A}},\mu_{\widetilde{B}})
  =\int_0^1\left(\int_X\left(1_{\mu_{\widetilde{A}}^{-1}(\alpha,1]\Delta\mu_{\widetilde{B}}^{-1}(\alpha,1]}(x)\right)d
    \sigma(x)\right)d\eta(\alpha).
\end{equation}

\begin{theorem} Assuming the integral well defined, \eqref{DMM} defines a metric, and
$\mu_{\widetilde{C}}$ is pointwise  between $\mu_{\widetilde{A}}$ and $\mu_{\widetilde{B}}$ if and only if the equality holds in the triangle inequality for this metric.
\label{thm7-3}
\end{theorem}

\begin{proof}
  By  Proposition \ref{propoDm} we have that for every $\alpha\in[0,1]$ the formula
\[
\int_X1_{A\Delta B}(x)d\sigma(x)
\]
defines a metric on $\mathcal P(X)$.
Thus \eqref{DMM} is an integral of metrics, and hence a metric.\smallskip

To prove the claim in the theorem we go along the lines of Proposition \ref{ineq-sigma-eta}, using \eqref{equa-q}
with $A$ replaced by $\mu_{\widetilde{A}}^{-1}(\alpha,1]$ and similarly for $B$ and $C$.  We can write:
  \[
    \begin{split}
 0& \le          D(\mu_{\widetilde{A}},\mu_{\widetilde{C}})+D(\mu_{\widetilde{C}},\mu_{\widetilde{B}})- D(\mu_{\widetilde{A}},\mu_{\widetilde{B}})=\\
&\hspace{-1cm}
=    2\int_0^1\left(\int_X\left(1_{\mu_{\widetilde{A}}^{-1}(\alpha,1]}(x)-1_{\mu_{\widetilde{C}}^{-1}(\alpha,1]}(x)\right)\left(1_{\mu_{\widetilde{B}}^{-1}(\alpha,1]}(x)-1_{\mu_{\widetilde{C}}^{-1}(\alpha,1]}(x)\right)d\eta(x)\right)d\sigma(\alpha)\\
&\hspace{-1cm}
=2\int_0^1\left(\left(\int_{\mu_{\widetilde{A}}^{-1}(\alpha, 1]\cap \mu_{\widetilde{B}}^{-1}(\alpha,1]}  \left(1-1_{\mu_{\widetilde{C}}^{-1}(\alpha,1]}(x)\right)d\eta(x)+\right.\right.\\
&     \left.\left. +
\int_{X\setminus(\mu_{\widetilde{A}}^{-1}(\alpha, 1]\cup \mu_{\widetilde{B}}^{-1}(\alpha,1])}  1_{\mu_{\widetilde{C}}^{-1}(\alpha,1]}(x)d\eta(x)\right)
        \right)d\sigma(\alpha) 
    \end{split}
  \]
  By the assumed properties on $d\sigma(\alpha)$ we have therefore equality in the triangle inequality if and only if
  \[
            \int_{\mu_{\widetilde{A}}^{-1}(\alpha, 1]\cap \mu_{\widetilde{B}}^{-1}(\alpha,1]}  \left(1-1_{\mu_{\widetilde{C}}^{-1}(\alpha,1]}(x)\right)d\eta(x)=
\int_{X\setminus(\mu_{\widetilde{A}}^{-1}(\alpha, 1]\cup \mu_{\widetilde{B}}^{-1}(\alpha,1])}  1_{\mu_{\widetilde{C}}^{-1}(\alpha,1]}(x)d\eta(x)=0
\]
and the end of the proof is as in the proof of Proposition \ref{ineq-sigma-eta}.
\end{proof}

\section{The hyperbolic numbers}
\setcounter{equation}{0}
\label{section2}
Complex numbers can be constructed as matrices of the form $\begin{pmatrix}a&-b\\ b&a\end{pmatrix}$ where $a,b\in\mathbb R$, and can be viewed (when $(a,b)\not=(0,0)$) as composition of an homothety and a rotation in the plane:
\[
\rho\begin{pmatrix}\cos \theta&-\sin\theta\\ \sin\theta&\cos \theta\end{pmatrix}
\]
Hyperbolic number in turn are symmetric matrices of the form

\begin{equation}
  \label{cond1}
  \begin{pmatrix}a&b\\b&a\end{pmatrix}=  \frac{1}{\sqrt{2}}\begin{pmatrix}1&1\\1&-1\end{pmatrix}\begin{pmatrix}a+b&0\\0&a-b\end{pmatrix}\frac{1}{\sqrt{2}}\begin{pmatrix}1&1\\1&-1\end{pmatrix}
\end{equation}
and can be seen when $a^2-b^2\not=0$ as composition of an homothety and an hyperbolic rotation
\[
\rho\begin{pmatrix}\cosh \theta&\sinh\theta\\ \sinh\theta&\cosh \theta\end{pmatrix}.
\]
Thus hyperbolic numbers from a family of pairwise commuting matrices; we refer to  \cite{MR3309382,MR3410909} for more information on these numbers.\smallskip

It will be convenient to set
\begin{equation}
  U=\frac{1}{\sqrt{2}}\begin{pmatrix}1&1\\1&-1\end{pmatrix}.
\end{equation}
Note that
\[
U=U^t\quad{\rm and}\quad U^2=I_2.
\]
We have
\begin{equation}
  \label{cond1111}
  \begin{pmatrix}a&-b\\-b&a\end{pmatrix}=  \frac{1}{\sqrt{2}}\begin{pmatrix}1&1\\1&-1\end{pmatrix}\begin{pmatrix}a-b&0\\0&a+b\end{pmatrix}\frac{1}{\sqrt{2}}\begin{pmatrix}1&1\\1&-1\end{pmatrix}
\end{equation}
and
\begin{equation}
  \begin{pmatrix}a&b\\b&a\end{pmatrix}  \begin{pmatrix}a&-b\\-b&a\end{pmatrix}=(a^2-b^2)I_2
  \end{equation}

Not every non-zero hyperbolic number is invertible but the formula
\begin{equation}
  \label{inverse}
\begin{pmatrix}a&b\\b&a\end{pmatrix}^{-1}=\frac{1}{a^2-b^2}\begin{pmatrix}a&-b\\-b&a\end{pmatrix},\quad a^2-b^2\not=0,
\end{equation}
shows in particular that the set of hyperbolic numbers for which $a^2-b^2\not=0$ form a multiplicative Abelian group of
matrices, a subgroup of which consists of the matrices for which $a^2-b^2=1$.\\

Thus:
\begin{lemma}
  The hyperbolic number satisfies
  \begin{equation}
    0\le   \begin{pmatrix}a&b\\b&a\end{pmatrix}\le \begin{pmatrix}1&0\\ 0&1\end{pmatrix}
  \end{equation}
  if and only if it holds that
  \begin{equation}
    \label{cond2}
    0\le a+b\le 1\quad and\quad 0\le a-b\le 1.
    \end{equation}
  \end{lemma}

  \begin{proof}
This is a direct consequence of \eqref{cond1}.
    \end{proof}

    Note that in the $(a,b)$ plane the set \eqref{cond2} is the square with vertices
\[
    (0,0), \,\,(1/2,1/2),\,\,(-1/2,-1/2)\quad{\rm and}\quad(1,0).
\]
    \begin{definition}
      \label{didi}
We denote by $\mathbb D$ the set of hyperbolic numbers satisfying \eqref{cond2}.
      \end{definition}
      \begin{figure}[h]
        \includegraphics[width=0.8\linewidth]{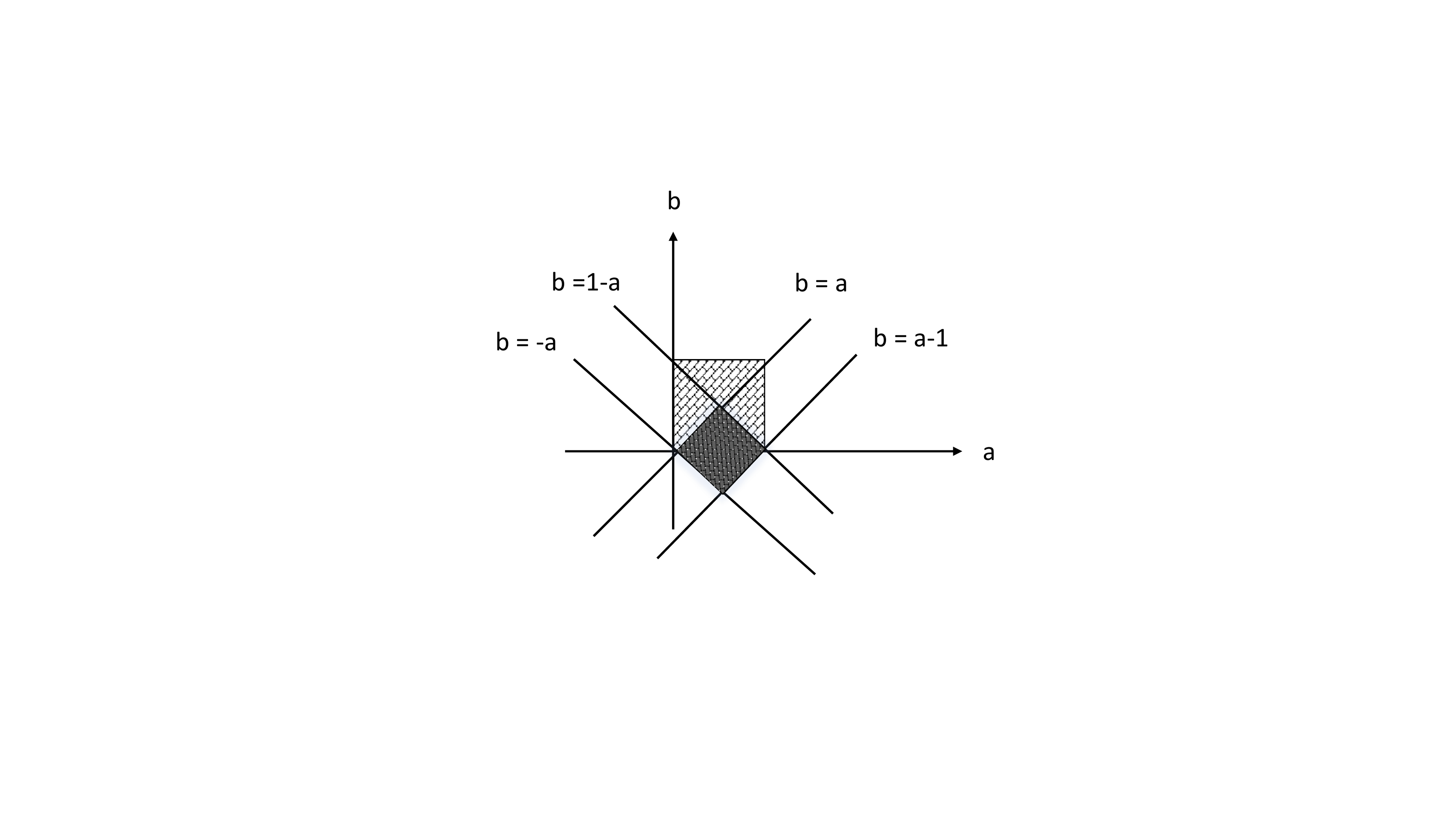}
        \vspace{-2cm}
\caption{The set $\mathbb D$}
\end{figure}
  For $z=\begin{pmatrix}a&b\\b&a\end{pmatrix}$ and $w=\begin{pmatrix}c&d\\d&c\end{pmatrix}$ in $\mathbb H$ we define
  \begin{equation}
    \label{new-york}
    \begin{split}
      z\vee w&=\frac{1}{\sqrt{2}}\begin{pmatrix}1&1\\1&-1\end{pmatrix}\begin{pmatrix}(a+b)\vee (c+d)
        &0\\0&(a-b)\vee (c-d)\end{pmatrix}\frac{1}{\sqrt{2}}\begin{pmatrix}1&1\\1&-1\end{pmatrix}\\
      &=\frac{1}{2}
      \begin{pmatrix}(a+b)\vee (c+d)+(a-b)\vee (c-d)        &(a+b)\vee (c+d)-(a-b)\vee (c-d)\\
(a+b)\vee (c+d)-(a-b)\vee (c-d)        &(a+b)\vee (c+d)+(a-b)\vee (c-d)
      \end{pmatrix}
      \end{split}
    \end{equation}
    and
      \begin{equation}
    \label{new-york}
    \begin{split}
      z\wedge w&=\frac{1}{\sqrt{2}}\begin{pmatrix}1&1\\1&-1\end{pmatrix}\begin{pmatrix}(a+b)\wedge (c+d)
        &0\\0&(a-b)\wedge (c-d)\end{pmatrix}\frac{1}{\sqrt{2}}\begin{pmatrix}1&1\\1&-1\end{pmatrix}\\
      &=\frac{1}{2}
      \begin{pmatrix}(a+b)\wedge (c+d)+(a-b)\wedge (c-d)        &(a+b)\wedge (c+d)-(a-b)\wedge (c-d)\\
(a+b)\wedge (c+d)-(a-b)\wedge (c-d)        &(a+b)\wedge (c+d)+(a-b)\wedge (c-d)
      \end{pmatrix}.
      \end{split}
    \end{equation}

    \begin{proposition}
The set $\mathbb H$ with the above functions $\vee$ and $\wedge$ is a lattice when endowed with the partial order of matrices.
      \end{proposition}
\begin{proof}
    It holds that
    \begin{equation}
z\wedge w\le     z\le z\vee w\quad{\rm and}\quad z\wedge w\le w\le z\vee w.
    \end{equation}
    We now discuss the uniqueness of the functions $\vee$ and $\wedge$. Given $z,w\in\mathbb D$ we consider positive hyperbolic numbers $m$ and $M$
    such that
    \begin{equation}
      \label{san-fransisco}
      m\le z\le M\quad{\rm and}\quad m\le w\le M
    \end{equation}
    We note that $m$ and $M$ are not unique, and two hyperbolic numbers $m_1$ and $m_2$ satisfying \eqref{san-fransisco} need not be comparable.
    But any $m$ and $M$ which satisfy \eqref{san-fransisco} will also satisfy
      \begin{equation}
        \label{kyoto}
        m\le  z\wedge w\quad{\rm and}\quad z\vee w\le M.
        \end{equation}

    \end{proof}

    As a corollary:

    \begin{corollary}
      In the above notation, $z\wedge w$ and $z\vee w$ are uniquely determined  to be respectively the largest and smallest hyperbolic numbers satisfying \eqref{kyoto}.
      \end{corollary}

      We now define the counterpart of an interval in the hyperbolic setting. Given two elements $z,w\in\mathbb H$, the characterization via \eqref{gamma-t} is not the
      one to consider here. Indeed the set
\[
\left\{c(t)=z\wedge v +t(z\vee w-z\wedge v),\,\, t\in[0,1]\right\}
\]
need not contain $z$ or $w$, as illustrated by the following example.
Take      
\begin{equation}
  \label{zw45}
z=U\begin{pmatrix}1&0\\0&2\end{pmatrix}U,\quad w=U\begin{pmatrix}3&0\\0&0\end{pmatrix}U.
\end{equation}
Then,
\begin{equation}
  \label{zw78}
z\wedge w=U\begin{pmatrix}1&0\\0&0\end{pmatrix}U,\quad z\vee w=U\begin{pmatrix}3&0\\0&2\end{pmatrix}U.
\end{equation}
Then, the interval
\[
[z\wedge w, z\vee w]=\left\{U\begin{pmatrix}1+2t&0\\0&2t\end{pmatrix}U,\,\, t\in[0,1]\right\}
\]
does not contain $z$ or $w$.\smallskip

Recall now that $\mathbb D$ was defined by condition \eqref{cond2} and denotes the set of positive
hyperbolic numbers less or equal to $I_2$.
      \begin{definition}
        Let $z,w\in\mathbb H$. We define the interval
        \begin{equation}
[z,w]_H=\left\{ c(\tau)=z\wedge w+\tau\left(z\vee w-z\wedge w\right),\,\, \tau\in\mathbb D\right\}.
          \end{equation}
        \label{london}
        \end{definition}

        \begin{proposition}
          $[z,w]_H$ can be characterized as:
\[
[z,w]_H=\left\{ v\in H\,\,:\,\, z\wedge w\le v\le z\vee w \right\}                
\]
\end{proposition}

\begin{proof}
  Let
  \[
\tau=U\begin{pmatrix}t_1&0\\0&t_2\end{pmatrix}U,\quad t_1,t_2\in[0,1]
    \]
  \[
    z=U\begin{pmatrix}\lambda_1&0\\0&\lambda_2\end{pmatrix}U,\quad
    w=U\begin{pmatrix}\mu_1&0\\0&\mu_2\end{pmatrix}U,\quad  {\rm and}\quad   v=U\begin{pmatrix}v_1&0\\0&v_2\end{pmatrix}U.
  \]
  We can write
  \begin{equation}
    \begin{split}
      c(\tau)&=U\begin{pmatrix}\underbrace{\lambda_1\wedge\mu_1+t_1(\lambda_1\vee\mu_1-\lambda_1\wedge\mu_1)}_{v_1(t_1)}&0\\
        0&\underbrace{\lambda_2\wedge\mu_2+t_2(\lambda_2\vee\mu_2-\lambda_2\wedge\mu_2)}_{v_2(t_2)}\end{pmatrix}U.
    \end{split}
    \label{reims}
    \end{equation}
    But, for $j=1,2$ and as $t_j$ varies from $0$ to $1$ we have that $v_j(t_j)$ varies from $\lambda_j\wedge \mu_j$ to $\lambda_j\vee \mu_j$.
    Hence, the representation  \eqref{reims} for $c(\tau)$ is equivalent to
    \[
 z\wedge w\le c(\tau)\le z\vee w .             
      \]
    
\end{proof}


                      \begin{definition}
                        Let $z,w,u\in \mathbb H$. We say that $u$ is between $z$ and $w$ if $u\in[z,w]_{\mathbb H}$, that is
                        \[
                          z\wedge w\le u\le z\vee w .
                          \]
                        \end{definition}

We note that the notion if betweenness is not transitive. Restle already had examples of lack of transitivity for sets. 

\begin{example}
  Take $z$ and $w$ as in \eqref{zw45}, with minimum and maximum as in \eqref{zw78} and
  \[
        u=U\begin{pmatrix}2&0\\0&0\end{pmatrix}U,\quad v=U\begin{pmatrix}5&0\\0&1/2\end{pmatrix}U.
    \]
    Then
\[
  u\wedge v=U\begin{pmatrix}2&0\\0&0\end{pmatrix}U,\quad u\vee v=U\begin{pmatrix}5&0\\0&1/2\end{pmatrix}U,
\]
and
\[
 z\wedge v=U\begin{pmatrix}1&0\\0&1/2\end{pmatrix}U,\quad z\vee v=U\begin{pmatrix}5&0\\0&2\end{pmatrix}U
\]
Then $u\in[z,w]_H$, $w\in[u,v]_H$ but $u\not \in[z,v]_H$.
\end{example}

 To conclude this section we note that we can write
    
    \begin{equation}
z=(a+b)P+(a-b)Q
\label{toronto}
\end{equation}
      where $P$ and $Q$ denote the orthogonal projections
      \begin{equation}
        \label{pqpq}
      P=\frac{1}{2}\begin{pmatrix}1\\1\end{pmatrix}\begin{pmatrix}1&1\end{pmatrix}\quad {\rm and}\quad
      Q=\frac{1}{2}\begin{pmatrix}1\\-1\end{pmatrix}\begin{pmatrix}1&-1\end{pmatrix}
      \end{equation}
      We further note that
      \begin{equation}
        P^2=P,\quad Q^2=Q,
        \end{equation}
      \begin{equation}
        \label{voltaire}
        PQ=QP=0.
      \end{equation}
      and
      \begin{equation}
        P+Q=I_2,
        \end{equation}
        and that
        \begin{equation}
          P=\frac{1}{2}\left(\begin{pmatrix}1&0\\0&1\end{pmatrix}+\begin{pmatrix}0&1\\1&0\end{pmatrix}\right)
        \end{equation}
        and
                \begin{equation}
          Q=\frac{1}{2}\left(\begin{pmatrix}1&0\\0&1\end{pmatrix}-\begin{pmatrix}0&1\\1&0\end{pmatrix}\right)
        \end{equation}
Representation \eqref{toronto} is called the idempotent representation of the hyperbolic number.
    In this work we chose to write hyperbolic numbers as matrices; one could also use the more traditional notation
    \[
      z=a+b{\mathbf k}
    \]
    where $\mathbf k\not\in\mathbb R$ satisfies ${\mathbf k}^2=1$ (in the matrix notation, we have ${\mathbf k}=\begin{pmatrix}0&1\\1&0\end{pmatrix}$).
    \section{$\mathbb D$-valued membership functions}
\setcounter{equation}{0}

\begin{definition}
  Let $X$ be a set. An hyperbolic-valued membership function on $X$ is a $\mathbb D$-valued map, i.e. a $\mathbb H$-valued map, say $M$, satisfying
  \begin{equation}
    0\le M(x)\le I_2,\quad x\in X.
    \end{equation}
  \end{definition}
\begin{theorem}
$M(x)$ is an hyperbolic-valued membership function if and only if there exist
  two membership functions  $\mu_{\widetilde{A_1}}$ and $\mu_{\widetilde{A_2}}$ 
corresponding to the fuzzy sets $\widetilde{A_1}$ and $\widetilde{A_2}$ respectively such that
\begin{equation}
  \label{totoche}
   M(x)=\frac{1}{2}\begin{pmatrix}\mu_{\widetilde{A_1}}(x)+\mu_{\widetilde{A_2}}(x)&\mu_{\widetilde{A_1}}(x)-\mu_{\widetilde{A_2}}(x)\\
    \mu_{\widetilde{A_1}}(x)-\mu_{\widetilde{A_2}}(x)&\mu_{\widetilde{A_1}}(x)+\mu_{\widetilde{A_2}}(x)\end{pmatrix}
  \end{equation}
\end{theorem}

\begin{proof}
 Following \eqref{cond1} we write
  \begin{equation}
    M(x)=\begin{pmatrix}a(x)&b(x)\\b(x)&a(x)\end{pmatrix}=  \frac{1}{\sqrt{2}}\begin{pmatrix}1&1\\1&-1\end{pmatrix}\begin{pmatrix}
      a(x)+b(x)&0\\0&a(x)-b(x)\end{pmatrix}\frac{1}{\sqrt{2}}\begin{pmatrix}1&1\\1&-1\end{pmatrix}.
  \end{equation}
  
    By \eqref{cond2} both $ a(x)+b(x)$ and $a(x)-b(x)$ take values in $[0,1]$ and hence are (classical) membership functions, corresponding to fuzzy sets
say  $\widetilde{A_1}$ and $\widetilde{A_2}$:
\[
  a(x)+b(x)=\mu_{\widetilde{A_1}}(x)\quad {\rm and}\quad a(x)-b(x)=\mu_{\widetilde{A_2}}(x).
\]
Thus
\begin{equation}
  \label{ma1a2}
M(x)=\frac{1}{\sqrt{2}}\begin{pmatrix}1&1\\1&-1\end{pmatrix}\begin{pmatrix}
     \mu_{\widetilde{A_1}}(x)&0\\0&\mu_{\widetilde{A_2}}(x)\end{pmatrix}\frac{1}{\sqrt{2}}\begin{pmatrix}1&1\\1&-1\end{pmatrix}.
  \end{equation}
Formula \eqref{totoche} follows.
  \end{proof}

We will use the notation
\[
  M(x)=    M_{\widetilde{A_1},\widetilde{A_2}}(x)
\]
and denote the fuzzy set as the pair $(\widetilde{A_1},\widetilde{A_2})$. One has

  \begin{equation}
    \label{lille}
    (\mu_{\widetilde{A_1}}\wedge\mu_{\widetilde{A_2}}) I_2\le 
M_{\widetilde{A_1},\widetilde{A_2}}(x)\le 
    (\mu_{\widetilde{A_1}}\vee\mu_{\widetilde{A_2}}) I_2
\end{equation}

It follows from \eqref{lille} that $M_{\widetilde{A_1},\widetilde{A_2}}$ defines a set ``between'' the intersection and the union of the two fuzzy sets
$\mu_{\widetilde{A_1}}\wedge\mu_{\widetilde{A_2}}$ and $\mu_{\widetilde{A_1}}\vee\mu_{\widetilde{A_2}}$.\smallskip  
  
We also note that we can rewrite  $M_{\widetilde{A_1},\widetilde{A_2}}(x)$ as the idempotent representation
\begin{equation}
  \label{st-ambroise}
  M_{\widetilde{A_1},\widetilde{A_2}}(x)=\mu_{\widetilde{A_1}}(x)P+\mu_{\widetilde{A_2}}(x)Q
    \end{equation}
where $P$ and $Q$ are as in \eqref{pqpq}.

          \begin{definition}
          Let $\alpha=\begin{pmatrix}\alpha_1&\alpha_2\\ \alpha_2&\alpha_1\end{pmatrix}$ be a positive hyperbolic number  less or equal to $I_2$. We define the $\alpha$-cut of the hyperbolic fuzzy set $M_{\widetilde{A_1},\widetilde{A_2}}$ to be
          \begin{equation}
            \label{haidak}
            \left\{x\in X\,\,;  M_{\widetilde{A_1},\widetilde{A_2}}(x)> \alpha\right\}
            \end{equation}
        \end{definition}

 By \eqref{cond1} we see  that  \eqref{haidak} is equivalent to
        \begin{eqnarray}
          \mu_{\widetilde{A_1}}(x)&>&\alpha_1+\alpha_2\\
          \mu_{\widetilde{A_2}}(x)&>&\alpha_1-\alpha_2,
                                      \end{eqnarray}
        corresponding to the (possibly empty) $\alpha$-cuts $(\widetilde{A_1})_{\alpha_1+\alpha_2}$ and $\alpha$-cuts $(\widetilde{A_2})_{\alpha_1-\alpha_2}$.
\begin{remark}
When $\mu_{\widetilde{A_1}}=1_{A_1}$ and $\mu_{\widetilde{A_2}}=1_{A_2}$ for some subsets $A_1$ and $A_2$ of $X$ we have
    \[
M_{\widetilde{A_1},\widetilde{A_2}}(x)=1_{A_1}(x)P+1_{A_2}(x)Q
\]
and
    \[
1_{A_1\cap A_2}I_2\le M_{\widetilde{A_1},\widetilde{A_2}}(x)\le 1_{A_1\cup A_2}I_2
      \]
    \end{remark}

    When
\[
\mu_{\widetilde{A_1}}(x)\ge \mu_{\widetilde{A_2}}(x),\quad \forall x\in X,
\]
both the functions
\begin{eqnarray}
  \mu_{ATA}(x)&=&\frac{\mu_{\widetilde{A_1}}(x)+\mu_{\widetilde{A_2}}(x)}{2}\\
            \nu_{ATA}(x)&=&\frac{\mu_{\widetilde{A_1}}(x)-\mu_{\widetilde{A_2}}(x)}{2}
\end{eqnarray}
are membership functions, and such that
\[
  \mu_{ATA}(x)+\nu_{ATA}(x)\le 1,
\]
and define an Atanassov intuitionistic fuzzy set.

\begin{definition}
  The hyperbolic fuzzy set is called an Atanassov hyperbolic fuzzy set if
  \[
    \mu_{\widetilde{A_1}}(x)\ge \mu_{\widetilde{A_2}}(x),\quad\forall x\in X
    \]
  \end{definition}

  \begin{proposition}
    The product of two Atanassov hyperbolic fuzzy sets is an
    Atanassov hyperbolic fuzzy set.
  \end{proposition}

  \begin{proof}
    Let $M_{\widetilde{A_1},\widetilde{A_2}}(x)$   and $M_{\widetilde{B_1},\widetilde{B_2}}(x)$ the two  Atanassov hyperbolic fuzzy set.
    Then
    \[
      \mu_{\widetilde{A_1}}(x)\ge \mu_{\widetilde{A_2}}(x),\quad {\rm and} \quad   \mu_{\widetilde{B_1}}(x)\ge \mu_{\widetilde{B_2}}(x),\quad\forall x\in X
            \]
            so that
            \[
              \mu_{\widetilde{A_1}}(x)\mu_{\widetilde{B_1}}(x)                \ge \mu_{\widetilde{A_2}}(x)\mu_{\widetilde{B_2}}(x),\quad\forall x\in X,
                    \]
                    and hence the answer.
    \end{proof}

\section{Properties of hyperbolic membership functions}
      \setcounter{equation}{0}
      In this section we consider the counterparts in the hyperbolic setting of the classical operators on fuzzy sets. We define

          \begin{equation}
M_{C(\widetilde{A_1},\widetilde{A_2})}(x)=I_2-M_{\widetilde{A_1},\widetilde{A_2}}(x)
            \end{equation}

and it is easy to verify that

        \begin{equation}
          M_{\widetilde{A_1},cl(\widetilde{A_1})}(x)=\frac{1}{2}\begin{pmatrix} 1&2\mu_{\widetilde{A_1}}(x)-1\\
           2 \mu_{\widetilde{A_1}}(x)-1&1\end{pmatrix}
          \end{equation}

 Furthermore, using \eqref{st-ambroise} and \eqref{voltaire} we have:

      \begin{proposition} 
Let $\widetilde{A_1},\widetilde{A_2},\widetilde{B_1}$ and $\widetilde{B_2}$ be fuzzy sets with membership functions 
$\widetilde{A_1},\widetilde{A_2},\widetilde{B_1}$ and $\widetilde{B_2}$ respectively. We have:
\begin{equation}
  M_{\widetilde{A_1},\widetilde{A_2}}(x)  M_{\widetilde{B_1},\widetilde{B_2}}(x)=\mu_{\widetilde{A_1}}(x)\mu_{\widetilde{B_1}}(x)P+\mu_{\widetilde{A_2}}(x)\mu_{\widetilde{B_2}}(x)Q.
  \end{equation}
\end{proposition}
Thus the matrix product of the hyperbolic membership functions $  M_{\widetilde{A_1},\widetilde{A_2}}$ and $M_{\widetilde{B_1},\widetilde{B_2}}$ corresponds to the algebraic product (see Section \ref{secfuzz} and \cite[\S 3.3. p. 33]{fuzz}) of the fuzzy sets
${\widetilde{A_1}}$ and ${\widetilde{B_1}}$ along $P$ and $\widetilde{A_2}$ and $\widetilde{B_2}$ along $Q$.\\

\begin{proposition}
  In the above notations it holds that:
  \begin{equation}
    M_{\widetilde{A_1},\widetilde{A_2}}(x)  M_{\widetilde{A_2},\widetilde{A_1}}(x)=\mu_{\widetilde{A_1}}(x)\mu_{\widetilde{A_2}}(x)I_2
    \end{equation}
    \end{proposition}
    \begin{proof}
      \[
        M_{\widetilde{A_1},\widetilde{A_2}}(x)  M_{\widetilde{A_2},\widetilde{A_1}}(x)=\mu_{\widetilde{A_1}}(x)\mu_{\widetilde{A_2}}(x)P+\mu_{\widetilde{A_2}}(x)\mu_{\widetilde{A_1}}(x)Q=\mu_{\widetilde{A_1}}(x)\mu_{\widetilde{A_2}}(x)I_2
        \]
      \end{proof}

By \eqref{new-york} we have:
\begin{equation}
  M_{\widetilde{A_1},\widetilde{A_2}}(x) \vee M_{\widetilde{B_1},\widetilde{B_2}}(x)=(\mu_{\widetilde{A_1}}(x)\vee \mu_{\widetilde{B_1}}(x))P+
  (\mu_{\widetilde{A_2}}(x)\vee\mu_{\widetilde{B_2}}(x))Q.
  \end{equation}

\section{Betweenness for hyperbolic-valued membership functions}
\setcounter{equation}{0}
\label{sec-end}
The proofs of the results in this section are easily adapted from the proofs in the scalar case by considering the idempotent decomposition, as in previous arguments in the paper, and we will not write out the details.
\begin{definition}
  Let $M_{\widetilde{A_1},\widetilde{A_2}}$, $M_{\widetilde{B_1},\widetilde{B_2}}$ and $M_{\widetilde{C_1},\widetilde{C_2}}$ three hyperbolic-valued membership functions defined on the
  set $X$.   One says that $M_{\widetilde{C_1},\widetilde{C_2}}$ is pointwise between $M_{\widetilde{A_1},\widetilde{A_2}}$ and $M_{\widetilde{B_1},\widetilde{B_2}}$ if
  \begin{equation}
    M_{\widetilde{A_1},\widetilde{A_2}}(x)\wedge M_{\widetilde{B_1},\widetilde{B_2}}(x)\le M_{\widetilde{C_1},\widetilde{C_2}}(x)\le
    M_{\widetilde{A_1},\widetilde{A_2}}(x)\vee M_{\widetilde{B_1},\widetilde{B_2}}(x),\quad x\in X.
    \end{equation}
  \end{definition}

In the setting of $\mathbb D$-valued membership functions
Lemma \ref{lemma4-2} and Proposition \ref{prop6-2} become:
  \begin{proposition}
Let $M_{\widetilde{A_1},\widetilde{A_2}}, M_{\widetilde{B_1},\widetilde{B_2}}$ and $M_{\widetilde{C_1},\widetilde{C_2}}$
be $\mathbb D$-valued membership functions. Then, $M_{\widetilde{C_1},\widetilde{C_2}}$ is pointwise between
$M_{\widetilde{A_1},\widetilde{A_2}}$ and $ M_{\widetilde{B_1},\widetilde{B_2}}$ if and only if
  \begin{equation}
    \label{nocrsip1111}
    M_{\widetilde{C_1},\widetilde{C_2}}(x)=M_{\widetilde{A_1},\widetilde{A_2}}(x)\wedge M_{\widetilde{B_1},\widetilde{B_2}}(x)+M_{\widetilde{Z_1},\widetilde{Z_2}}(x)
  \end{equation}
  where $M_{\widetilde{C_1},\widetilde{C_2}}(x)$ is a $\mathbb D$-valued membership function satisfying
  \begin{equation}
    \label{nocrsip2222}
    M_{\widetilde{Z_1},\widetilde{z_2}}(x)\le M_{\widetilde{A_1},\vee{A_2}}(x)\vee M_{\widetilde{B_1},\widetilde{B_2}}(x)-M_{\widetilde{A_1},\widetilde{A_2}}(x)
    \wedge M_{\widetilde{B_1},\widetilde{B_2}}(x),\quad x\in X.
  \end{equation}
  \label{prop6-2322222}
  \end{proposition}

  Furthermore, the idempotent decomposition \eqref{toronto} gives:

  \begin{proposition}
In the notation of the previous proposition, $M_{\widetilde{C_1},\widetilde{C_2}}$ is pointwise between
$M_{\widetilde{A_1},\widetilde{A_2}}$ and $ M_{\widetilde{B_1},\widetilde{B_2}}$ if and only if
$\widetilde{C_1}$ and $\widetilde{C_2}$ are pointwise between
$\widetilde{A_1}$ and $\widetilde{B_1}$ and $\widetilde{A_2}$ and $\widetilde{B_2}$ respectively.
    \end{proposition}

  The hyperbolic counterparts of Definition \ref{orleans} and Theorem \ref{bretagne} in the hyperbolic setting are:
    
    \begin{definition}
      Let $\mathbf{a}\in\mathbb D$. The $\mathbf a$-cut associated to the $\mathbb D$-valued membership function $M_{\widetilde{A_1},\widetilde{A_2}}$ is the set of
      elements $M_{\widetilde{A_1},\widetilde{A_2}}^{-1}$.
      \end{definition}

      \begin{theorem}
$M_{\widetilde{C_1},\widetilde{C_2}}$ is ${\mathbf a}$-between $M_{\widetilde{A_1},\widetilde{A_2}}$ and $M_{\widetilde{B_1},\widetilde{B_2}}$ 
      if and only if $M_{\widetilde{C_1},\widetilde{C_2}}$ is pointwise between $M_{\widetilde{A_1},\widetilde{A_2}}$  and $M_{\widetilde{B_1},\widetilde{B_2}}$ 
        \end{theorem}

        We conclude with a counterpart of Theorem \ref{thm7-3} for hyperbolic-valued membership functions. The novelty is what one now needs the
        $\mathbb H$-valued counterpart of a distance to get a triangle equality. Here too the proof is easy, going via the idempotent decomposition \eqref{toronto},
        and will be omitted. With $D(\mu_{\widetilde{A}},\mu_{\widetilde{B}})$ as in \eqref{DMM} we define
        \begin{equation}
          D_{\mathbb H}(M_{\widetilde{A_1},\widetilde{A_2}},M_{\widetilde{B_1},\widetilde{B_2}} )=U\begin{pmatrix} 
          D(\mu_{\widetilde{A_1}},\mu_{\widetilde{B_1}}) &0\\         0&               D(\mu_{\widetilde{A_2}},\mu_{\widetilde{B_2}})\end{pmatrix}U.
      \end{equation}

      \begin{theorem} Assuming the integral well defined,
        $M_{\widetilde{C_1},\widetilde{C_2}}$ is pointwise between $M_{\widetilde{A_1},\widetilde{A_2}}$  and $M_{\widetilde{B_1},\widetilde{B_2}}$  if and only if the equality holds
        \[
          D_{\mathbb H}(M_{\widetilde{A_1},\widetilde{A_2}},M_{\widetilde{B_1},\widetilde{B_2}} )=D_{\mathbb H}(M_{\widetilde{A_1},\widetilde{A_2}},M_{\widetilde{C_1},\widetilde{C_2}} )+
          D_{\mathbb H}(M_{\widetilde{C_1},\widetilde{C_2}},M_{\widetilde{B_1},\widetilde{B_2}} ).
          \]
\label{thm7-4}
\end{theorem}
        
        \bibliographystyle{plain}

\begin{thebibliography}{10}

\bibitem{aizenberg2013multi}
I.~Aizenberg, N.~Aizenberg, and J.~Vandewalle.
\newblock {\em Multi-valued and universal binary neurons: Theory, learning and
  applications}.
\newblock Kluwer Academic Publishers, 2000.

\bibitem{MR2002b:47144}
D.~Alpay.
\newblock {\em The {S}chur algorithm, reproducing kernel spaces and system
  theory}.
\newblock American Mathematical Society, Providence, RI, 2001.
\newblock Translated from the 1998 French original by Stephen S. Wilson,
  Panoramas et Synth\`eses.

\bibitem{MR4537587}
D.~Alpay, K.~Diki, and M.~Vajiac.
\newblock A note on the complex and bicomplex valued neural networks.
\newblock {\em Appl. Math. Comput.}, 445:Paper No. 127864, 12, 2023.

\bibitem{MR4302453}
D.~Alpay and P.~Jorgensen.
\newblock New characterizations of reproducing kernel {H}ilbert spaces and
  applications to metric geometry.
\newblock {\em Opuscula Math.}, 41(3):283--300, 2021.

\bibitem{MR3309382}
D.~Alpay, M.~Luna-Elizarrar{\'a}s, M.~Shapiro, and D.C. Struppa.
\newblock {\em Basics of functional analysis with bicomplex scalars, and
  bicomplex {S}chur analysis}.
\newblock Springer Briefs in Mathematics. Springer, Cham, 2014.

\bibitem{MR3651492}
D.~Alpay, M.E. Luna-Elizarrar\'{a}s, and M.~Shapiro.
\newblock Kolmogorov's axioms for probabilities with values in hyperbolic
  numbers.
\newblock {\em Adv. Appl. Clifford Algebr.}, 27(2):913--929, 2017.

\bibitem{aron}
N.~Aronszajn.
\newblock Theory of reproducing kernels.
\newblock {\em Trans. Amer. Math. Soc.}, 68:337--404, 1950.

\bibitem{MR852871}
K.T. Atanassov.
\newblock Intuitionistic fuzzy sets.
\newblock {\em Fuzzy Sets and Systems}, 20(1):87--96, 1986.

\bibitem{chen-pham}
G.~Chen and T.~Pham.
\newblock {\em {Introduction to fuzzy sets, fuzzy logic and fuzzy control
  systems}}.
\newblock {CRC Press LLC}, 2001.

\bibitem{MR2522441}
Shihyen Chen, Bin Ma, and Kaizhong Zhang.
\newblock On the similarity metric and the distance metric.
\newblock {\em Theoret. Comput. Sci.}, 410(24-25):2365--2376, 2009.

\bibitem{MR224391}
J.~A. Goguen.
\newblock {$L$}-fuzzy sets.
\newblock {\em J. Math. Anal. Appl.}, 18:145--174, 1967.

\bibitem{hays1958approach}
W.L. Hays.
\newblock An approach to the study of trait implication and trait similarity.
\newblock In R.~Tagiuri and L.~Petrullo, editors, {\em Person perception and
  interpersonal behavior}, pages 289--299. Stanford University Press Stanford,
  1958.

\bibitem{MR0345715}
Y.~Horibe.
\newblock A note on entropy metrics.
\newblock {\em Information and Control}, 22:403--404, 1973.

\bibitem{jaccard1901distribution}
P.~Jaccard.
\newblock Distribution de la flore alpine dans le bassin des dranses et dans
  quelques r{\'e}gions voisines.
\newblock {\em Bull Soc Vaudoise Sci Nat}, 37:241--272, 1901.

\bibitem{kobayashi2012hyperbolic}
M.~Kobayashi.
\newblock Hyperbolic {H}opfield neural networks.
\newblock {\em IEEE transactions on neural networks and learning systems},
  24(2):335--341, 2013.

\bibitem{kuroe}
Y.~Koroe.
\newblock A model of complex-valued associative memories and its dynamics.
\newblock In Akira Hirose, editor, {\em Complex-valued neural networks}, pages
  57--79. World Scientific, 2003.

\bibitem{levandowsky1971distance}
M.~Levandowsky and D.~Winter.
\newblock Distance between sets.
\newblock {\em Nature}, 234(5323):34--35, 1971.

\bibitem{MR1356718}
H.X. Li and V.C. Yen.
\newblock {\em Fuzzy sets and fuzzy decision-making}.
\newblock CRC Press, Boca Raton, FL, 1995.

\bibitem{zbMATH03335866}
J.~{\L}ukasiewicz.
\newblock Selected works. {Edited} by {L}. {Borkowski}.
\newblock Studies in {Logic} and the {Foundations} of {Mathematics}.
  {Amsterdam} etc.: {North}-{Holland} {Publishing} {Company}; {Warszawa}: {PWN}
  - {Polish} {Scientific} {Publishers}. xii, 405 pages, 1970.

\bibitem{MR3410909}
M.E. Luna-Elizarrar\'{a}s, M.~Shapiro, D.C. Struppa, and A.~Vajiac.
\newblock {\em Bicomplex holomorphic functions}.
\newblock Frontiers in Mathematics. Birkh\"{a}user/Springer, Cham, 2015.
\newblock The algebra, geometry and analysis of bicomplex numbers.

\bibitem{marczewski1958certain}
E.~Marczewski and H.~Steinhaus.
\newblock On a certain distance of sets and the corresponding distance of
  functions.
\newblock In {\em Colloquium Mathematicum}, volume~6, pages 319--327. Instytut
  Matematyczny Polskiej Akademii Nauk, 1958.

\bibitem{zbMATH03194216}
E.~Marczewski and H.~Steinhaus.
\newblock On the systematik distance of biotopes.
\newblock {\em Zastosow. Mat.}, 4:195--202, 1959.

\bibitem{fuzz}
C.~Mohan.
\newblock {\em An introduction to fuzzy set theory and fuzzy logic}.
\newblock {MV Learning}, second edition, 2019.

\bibitem{zbMATH01202379}
W.~Pedrycz and F.~Gomide.
\newblock {\em An introduction to fuzzy sets: analysis and design. {With} a
  foreword by {Lotfi} {A}. {Zadeh}}.
\newblock Cambridge, MA: MIT Press, 1998.

\bibitem{pereverzyev2022introduction}
S.~Pereverzyev.
\newblock {\em An Introduction to Artificial Intelligence Based on Reproducing
  Kernel Hilbert Spaces}.
\newblock Springer Nature, 2022.

\bibitem{rajski_1961}
C.~Rajski.
\newblock A metric space of discrete probability distributions.
\newblock {\em {Information and Control}}, 4:371--377, 1961.

\bibitem{ralescu1984probability}
A.L. Ralescu and D.A. Ralescu.
\newblock Probability and fuzziness.
\newblock {\em Information Sciences}, 34(2):85--92, 1984.

\bibitem{restle1959metric}
F.~Restle.
\newblock A metric and an ordering on sets.
\newblock {\em Psychometrika}, 24(3):207--220, 1959.

\bibitem{saitoh}
S.~Saitoh.
\newblock {\em Theory of reproducing kernels and its applications}, volume 189.
\newblock Longman scientific and technical, 1988.

\bibitem{MR10:133e}
C.E. Shannon.
\newblock A mathematical theory of communication.
\newblock {\em Bell System Tech. J.}, 27:379--423, 623--656, 1948.

\bibitem{MR632835}
T.~A. Springer.
\newblock {\em Linear algebraic groups}, volume~9 of {\em Progress in
  Mathematics}.
\newblock Birkh\"auser Boston, Mass., 1981.

\bibitem{tversky1977features}
A.~Tversky.
\newblock Features of similarity.
\newblock {\em Psychological review}, 84(4):327, 1977.

\bibitem{zbMATH03226241}
L.A. Zadeh.
\newblock Fuzzy sets.
\newblock {\em Inf. Control}, 8:338--353, 1965.

\bibitem{MR1174743}
H.-J. Zimmermann.
\newblock {\em Fuzzy set theory---and its applications}.
\newblock Kluwer Academic Publishers, Boston, MA, second edition, 1992.
\newblock With a foreword by L. A. Zadeh.

\bibitem{zwick1987measures}
R.~Zwick, E.~Carlstein, and D.V. Budescu.
\newblock Measures of similarity among fuzzy concepts: A comparative analysis.
\newblock {\em International journal of approximate reasoning}, 1(2):221--242,
  1987.

\end{thebibliography}
\def\cprime{$'$} \def\cprime{$'$} \def\cprime{$'$}
  \def\lfhook#1{\setbox0=\hbox{#1}{\ooalign{\hidewidth
  \lower1.5ex\hbox{'}\hidewidth\crcr\unhbox0}}} \def\cprime{$'$}
  \def\cprime{$'$} \def\cprime{$'$} \def\cprime{$'$} \def\cprime{$'$}
  \def\cprime{$'$}

  \end{document}